\documentclass[11pt]{article}
\font\smallit=cmti10

\usepackage[english]{babel}
\usepackage{amssymb,amsthm,amsmath,amsfonts,mathrsfs}
\usepackage{pst-all,pst-3dplot,pstricks,pstricks-add,pst-math,pst-xkey}
\usepackage{graphicx}
\usepackage{xstring}
\usepackage{pgffor}
\usepackage{latexsym}
\usepackage{ifthen, bigstrut}
\usepackage{accents}
\usepackage{hyperref}
\usepackage{breakurl}
\usepackage{fixmath}
\usepackage{xstring}

\usepackage{bbold}
\usepackage{enumerate}
\usepackage{epsfig}
\usepackage{subfigure}
\usepackage{skak}
\usepackage{chessboard}
\usepackage{chessfss}
\usepackage{pst-all}
\usepackage{lscape}
\usepackage{courier}
\usepackage[all,cmtip]{xy}

\usepackage{wasysym}

\newtheorem{theorem}{Theorem}

\newtheorem{lemma}[theorem]{Lemma}

\tikzstyle{hackennode}=[draw,circle,inner sep=0,minimum size=2pt, fill=white]
\tikzstyle{hackenline}=[line width=1.5pt]
\tikzstyle{hacken2remove}=[hackenline,green!30!black]
\tikzstyle{hackenroot}=[hackenline,green!30!purple]

\newcommand{\fixture}[1]{%
\StrLeft{#1}{1}[\firstchar]%
\ifdefstring{\firstchar}{R}{%
\includegraphics[scale=0.2]{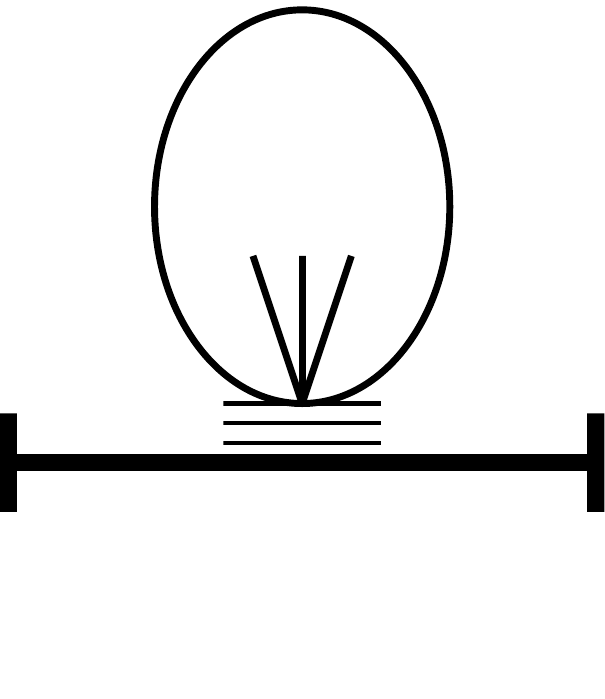}%
}{\ifdefstring{\firstchar}{B}{%
\includegraphics[scale=0.2]{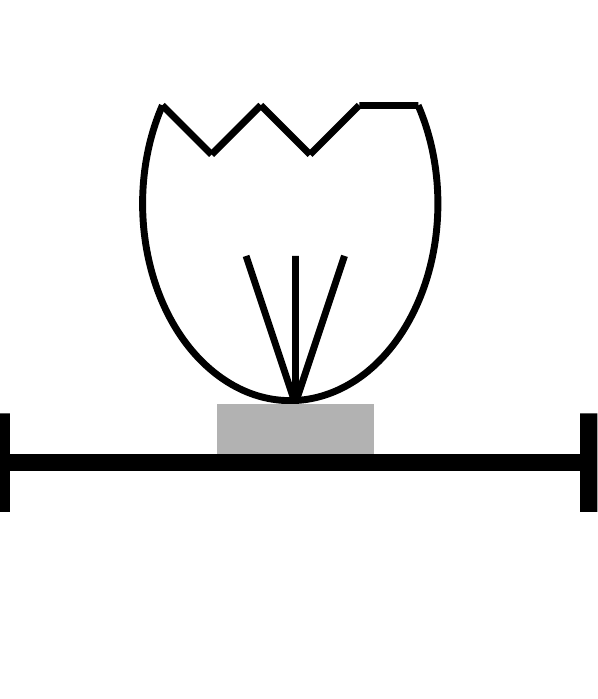}%
}{%
\includegraphics[scale=0.2]{\firstchar.pdf}%
}%
}%
\StrLen{#1}[\stringlength]%
\ifnum\stringlength>1%
\foreach \i in {2,...,\stringlength}{%
\StrChar{#1}{\i}[\currentchar]%
\includegraphics[scale=0.2]{\currentchar.pdf}%
}%
\fi%
}

\let\leftmoon\relax

\usepackage{mathabx}

\usepackage{bbold}
\newcommand{\bs}{\boldsymbol}

\renewcommand{\inf}{\ensuremath{{\bs\infty }}}

\newcommand{\Imi}{\mathbb{Im}^\inf}
\renewcommand{\Im}{\mathbb{nim}}
\newcommand{\imm}[1]{S_{#1}}
\newcommand{\prot}[1]{P_{#1}}

\newcommand{\moon}{\scalebox{1.1}{\ensuremath{\leftmoon}}}

\newcommand{\mex}{\ensuremath{\mathrm{mex}}}

\newcommand{\GL}{{G^\mathcal{L}}}

\theoremstyle{definition}

\newtheorem{observation}[theorem]{Observation}

\makeatletter
\newcommand{\oset}[3][0ex]{%
  \mathrel{\mathop{#3}\limits^{
    \vbox to#1{\kern-1.5\ex@
    \hbox{$\scriptstyle#2$}\vss}}}}
\makeatother

\emergencystretch=\maxdimen
\begin{document}
\begin{center}
\uppercase{\bf A complete solution for a nontrivial ruleset with entailing moves}
\vskip 6pt
{\bf Urban Larsson\footnote{larsson@iitb.ac.in}}\\
{\smallit Indian Institute of Technology Bombay, India}\\
{\bf Richard J.~Nowakowski\footnote{r.nowakowski@dal.ca}}\\
{\smallit Department of Mathematics and Statistics, Dalhousie University, Canada}\\
{\bf Carlos P. Santos\footnote{Under the scope of
UIDB/00297/2020 and UIDP/00297/2020; cmf.santos@fct.unl.pt}}\\
{\smallit Center for Mathematics and Applications (NovaMath), FCT NOVA, Portugal}\\
\end{center}

\begin{abstract}
Combinatorial Game Theory typically studies sequential rulesets with perfect information where two players alternate moves. There are rulesets with {\em entailing  moves} that break the alternating play axiom and/or restrict the other player's options within the disjunctive sum components. Although some examples have been analyzed in the classical work \emph{Winning Ways}, such rulesets usually fall outside the scope of the established normal play mathematical theory. At the first Combinatorial Games Workshop at MSRI, John H. Conway proposed that an effort should be made to devise some nontrivial ruleset with entailing moves that had a complete analysis. Recently, Larsson, Nowakowski, and Santos proposed a more general theory, {\em affine impartial}, which facilitates the mathematical analysis of impartial rulesets with entailing  moves. Here, by using this theory, we present a complete solution for a nontrivial ruleset with entailing moves.
\end{abstract}

\vspace{-0.5cm}
\section{Introduction}
The theory of disjunctive sums of combinatorial games was introduced by Conway~\cite{Con1976} and further expanded by Berlekamp, Conway and Guy in ``Winning Ways'' \cite{BerleCG1982}. The main point of the theory is that if a ruleset decomposes into components, then the analysis becomes easier. Each component is assigned a theoretical value, which is an abstract concept that is not tied to the ruleset. A position is a sum of individual components. An important fact is that the players move alternately in the position, but not necessarily in the components.

Winning Ways considers many types of rulesets which are not fully covered by this theory. Half of Chapter 12 involves impartial rulesets with \textit{entailing moves}. No theory is given. Some rulesets are considered, although none are solved. With an entailing move, if a certain condition occurs, the options of the next player are reduced -- for example, an entailing move may force the opponent to play on a certain pile. In this document, the game forms of entailing moves are expressed with help of the symbols $\infty$ (an unconditional Left win) and $\overline{\infty}$ (an unconditional Right win). Of special interest is the moon value, $\moon=\{\infty\,|\,\overline{\infty}\}$, where each player has a terminating move. This is part of a general theory \cite{LNS}, explained further on. A special note on terminology: since the moon is the only new affine value \cite{LNS} (adjoined to the nimbers), we designate the term {\em Grundy-value} for our generalized Grundy-value.

It was noted by the authors that entailing moves also occurred in \textsc{nimstring} (the impartial version of \textsc{dots and boxes}, see Chapter 16) and other rulesets. In those particular cases, the authors of Winning Ways used the designation \textit{complimenting moves}, where the players `carry on' the moves, keeping the turn to play. Given a game $G$, such a move has the form $G^L=\{\inf|G^{LR}\}$ or $G^R=\{G^{RL}|\overline{\infty}\}$. Once there is a lethal threat expressed by the infinity symbol, there is an automatic ``jump'' from $G$ to $G^{LR}$ or to $G^{RL}$. These moves can also be seen as moves that reduce options, since the player is forced to respond locally in a certain way to protect himself from an infinitely large threat. Hence, complementing moves are particular cases of entailing moves. Here, we refer to complimenting moves as \emph{carry-on} moves.

Although \textsc{nimstring} and \textsc{dots and boxes} have received attention \cite{Berle2000,Bremetal,Ish1,Ish2}, until recently little progress has been made towards a general theory. There are only two papers which mention entailing moves: in 1996, \cite{West1996}, which is a computer analysis of {\sc top entails} heaps, ranging in size up to $600,000$, and no regularities were discovered; and, in 2002, \cite{Elkies2002}, which considers pawn endgames in {\sc chess}. In the latter, entailing moves avoid losing immediately but no other theory is needed.

At the first Combinatorial Games Workshop at MSRI (1996), John H. Conway proposed that an effort should be made to devise some nontrivial ruleset with entailing moves that had a complete analysis.
In this paper, we introduce the \textsc{christmas lights' fixture}, which has carry-on moves, and we give the complete analysis in Section \ref{sec:Analysis}.

One reason why a complete analysis is possible is that, in \cite{LNS}, we show that impartial games, with entailing moves (and, in particular, carry-on moves), can be incorporated into one theory that extends impartial normal play structure. We review the \emph{affine impartial normal play} theory in Section \ref{sec:Review}.\\

\subsection{{\sc christmas lights' fixture}}
The ruleset \textsc{christmas lights' fixture} is inspired by the Christmas season. In a typical family home, the Christmas tree is decorated by strings of lights which we call ``fixtures''. With use, some parts of the fixtures tend to become damaged. These damaged parts may have either broken bulb sockets (broken bases), which must be removed, or broken-but-replaceable bulbs. Observe that a broken-but-replaceable bulb has a ``live'' socket in good shape. For example, in the fixture\\

\noindent
\scalebox{0.63}{\fixture{PSYRRBRCOGRBYSRRR}}

\vspace{-0.7cm}
\hspace{1.45cm}$\underbrace{\quad\quad\quad\quad\quad\quad\quad\,\,\,\,\,}_{\text{Damaged}}$\hspace{2.2cm}$\underbrace{\quad\quad\quad\quad\,\,}_{\text{Damaged}}$\hspace{1.5cm}$\underbrace{\quad\quad\quad\quad\quad\quad}_{\text{Damaged}},$\\

\noindent
a working bulb is colored, a replaceable bulb is transparent, and a broken socket is indicated by a broken bulb and a black socket.
This fixture has three disjoint damaged zones.\\

After Christmas, the family calls two electricians to repair a fixture. Of course, they want to replace all the replaceable bulbs and remove all the broken sockets. This is a quiet time, so the electricians play the following game in which they alternate moves. Since they only work on the damaged parts of the fixture, the whole fixture may be seen as a disjunctive sum where the disjoint components are the damaged parts, separated by working bulbs. The previous example corresponds to the disjunctive sum

\vspace{0.3cm}
\begin{center}\scalebox{0.63}{\fixture{RRBR}}\raisebox{0.38cm}{$\,\,\,\,+\,\,\,\,$}\scalebox{0.63}{\fixture{RB}}\raisebox{0.38cm}{$\,\,\,\,+\,\,\,\,$}\scalebox{0.63}{\fixture{RRR}}\end{center}

\begin{figure}[htb!]
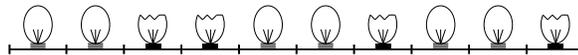

\caption{A damaged part is a component of a disjunctive sum.}
\label{fig:*1}
\begin{center}
\scalebox{0.63}{\fixture{RRBBRRBRRB}}
\end{center}
\end{figure}

\vspace{0.2cm}
\noindent
\emph{Positions}: A \textsc{christmas lights' fixture} position is a fixture with some damaged parts. Figure~\ref{fig:*1} shows a possible damaged part to be used to exemplify how the moves are made.

\noindent
\emph{Moves}: There are three types of moves:\\ 

\begin{figure}[htb!]
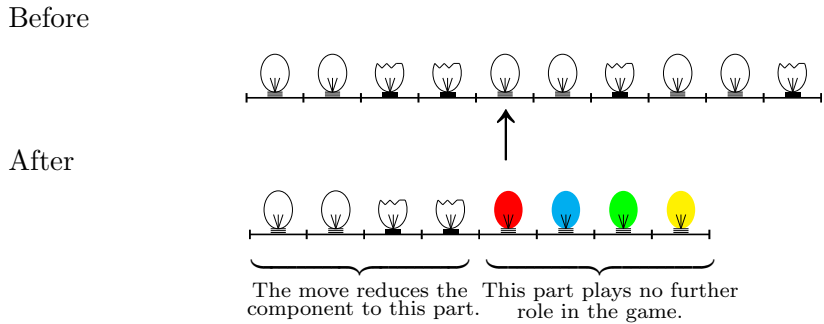

\caption{Replacing a replaceable bulb.}
\label{fig:*2}
Before
\begin{center}\scalebox{0.63}{\fixture{RRBBRRBRRB}}\end{center}

\vspace{-0.5cm}\hspace{6.42cm}\scalebox{2}{$\uparrow$}

\vspace{-0.2cm}
After\\

\vspace{-0.25cm}
\hspace{3.2cm}\scalebox{0.63}{\fixture{RRBBSCGY}}
\end{figure}

\noindent (1) If an electrician chooses a replaceable bulb of a component to play, then she fixes that bulb  and everything on that component to the right (away from the plug). That is, the electrician replaces replaceable bulbs and removes broken sockets to the right -- Figure~\ref{fig:*2}.\\

\vspace{-0.9cm}
\hspace{2.65cm}$\underbrace{\quad\quad\quad\quad\quad\quad\quad\,\,\,\,}_{\substack{\text{The move reduces the} \\ \text{ component to this part.}}}$\hspace{-0.25cm}$\underbrace{\quad\quad\quad\quad\quad\quad\quad\quad}_{\substack{\text{\;\;
This part plays no further} \\ \text{role in the game.}}}$\\

\begin{figure}[htb!]
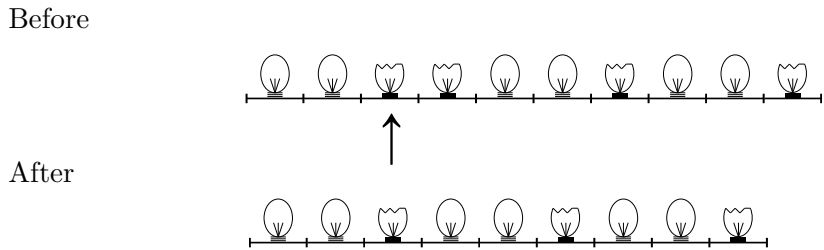

\caption{Removing a broken socket without getting a shock.}
\label{fig:*3}

Before
\begin{center}\scalebox{0.63}{\fixture{RRBBRRBRRB}}\end{center}

\vspace{-0.45cm}\hspace{4.9cm}\scalebox{2}{$\uparrow$}

\vspace{-0.1cm}
After\\

\vspace{-0.3cm}
\hspace{3.2cm}\scalebox{0.63}{\fixture{RRBRRBRRB}}
\end{figure}

\vspace{0.5cm}
\noindent (2) If an electrician chooses a broken socket to play, and if the adjacent bulbs are not both replaceable (with live sockets), then the electrician can remove it and connect the two ends together, passing the turn to the opponent -- Figure~\ref{fig:*3}.\\

\noindent (3) However, if an electrician chooses a broken socket to play, and if the adjacent bulbs are both broken-but-replaceable, then, when connecting the two ends, the electrician gets a mild electric shock. As a consequence of the shock, the electrician \emph{must move again}, which can be done on any component. This is the \textit{carry-on} rule. Removing a broken socket that is adjacent to only one or no broken-but-replaceable bulb does not trigger a shock nor a carry-on move -- Figure~\ref{fig:*4}.

\begin{figure}[htb!]
\caption{Removing a broken socket, and getting a mild shock.}
\label{fig:*4}

Before
\begin{center}\scalebox{0.63}{\fixture{RRBBRRBRRB}}\end{center}

\vspace{-0.5cm}\hspace{7.96cm}\scalebox{2}{$\uparrow$}
\vspace{0.2cm}

After (the player must play again on any component)\\

\vspace{0.3cm}
\hspace{3.2cm}\scalebox{0.63}{\fixture{RRBBRRRRB}}

\vspace{-1.3cm}\hspace{7.6cm}\scalebox{2.2}{\textcolor[rgb]{1.00,0.50,0.25}{\lightning}}\\
\end{figure}
\vspace{.3cm}
\noindent
\emph{Winning condition}: An electrician who no longer has broken-but-replaceable bulbs to replace or broken bulb sockets to remove loses the game (normal play convention).\\

\vspace{-0.3cm}
Moves of type (1) are similar to {\sc green hackenbush}. Indeed, there are even deeper connections.

\section{Review of theory, an exposition, and some notation}\label{sec:Review}
Because working bulbs take no part of the game, we will henceforth designate the terminology {\em socket} for a broken socket, and {\em bulb} for a broken-but-replaceable bulb.

The proofs in this paper will focus on concatenations with bulbs and sockets that are reminiscent of \emph{ordinal sums}. To motivate the concept, let us consider the classic ruleset \textsc{green hackenbush}, which is different from \textsc{nim} \cite{BerleCG1982}. The next left diagram shows a position where, after removing the edge with the label~\emph{a}, four more edges on the top disappear. In this example, it makes sense to consider the decomposition shown in the right diagram.

\begin{center}
\begin{tikzpicture}[scale=0.4]

    \draw [dotted] (-2.5,0) -- (2.5,0); 

    \node[hackennode] (root1)  at (0,0) {};
    \node[hackennode] (root2)  at (0, 1.2) {};
    \node[hackennode] (root3)  at (0, 2.4) {};
    \node[hackennode] (root4)  at (0, 3.6) {};
    \node[scale=1] at (0.4, 3.0) {$\,\,a$};
    \node[hackennode] (root5a)  at (-0.85, 4.45) {};
    \node[hackennode] (root5b)  at (0.85, 4.45) {};
    \node[hackennode] (root6)  at (0, 5.3) {};

    \draw[hackenline,green]
        (root1) -- (root2) -- (root3) -- (root4) -- (root5a) -- (root6)
        (root4) -- (root5b) -- (root6);

\end{tikzpicture}
\hspace{1.25cm}
\begin{tikzpicture}[scale=0.6]

    \node[scale=1] at (-2.7, 1.85) {$G=$};
    \node[hackennode] (root1a)  at (-2,0) {};
    \node[hackennode] (root2a)  at (-2, 1.2) {};
    \node[hackennode] (root3a)  at (-2, 2.4) {};
    \node[hackennode] (root4a)  at (-2, 3.6) {};

    \node[scale=1] at (0.5, 1.85) {$H=\,\,\,$};
    \node[hackennode] (root1b)  at (2, 1) {};
    \node[hackennode] (root2b)  at (1.15, 1.85) {};
    \node[hackennode] (root3b)  at (2.85, 1.85) {};
    \node[hackennode] (root4b)  at (2, 2.7) {};

    \draw[hackenline,green]
        (root1a) -- (root2a) -- (root3a) -- (root4a)
        (root1b) -- (root2b) -- (root4b)
        (root1b) -- (root3b) -- (root4b);

\end{tikzpicture}
\end{center}

If a player moves in $G$, then $H$ disappears; if a player moves in $H$, then nothing happens to $G$. The intuitive understanding in play is that $H$ is eliminated when a player moves in $G$. This idea leads to the concept of ordinal sum of two games $G:H$, where a player may move in either $G$ (base) or $H$ (subordinate), with the additional constraint that any move on $G$ completely annihilates the component $H$. The recursive definition is $$G:H=\left\{G^\mathcal{L},G:H^\mathcal{L}\,|\,G^\mathcal{R},G:H^\mathcal{R}\right\}.$$
It is crucial to remember the Colon Principle, i.e., if $H= H'$, then $G:H= G:H'$, but, in general, $H:G\ne H':G$. There are many works where the ordinal sum is the key concept (the classical \cite{AlberNW2007,BerleCG1982,Con1976,Siegel2013}, or the more recent \cite{CNS2018,CHNS2021,FNSW2015,FNS2018}).
Typically, one begins by proving that the value of $G:H$ can be determined, provided that $G$ and $H$ are minimal in some respect. Then, given a composition of ordinal sums
$$G=G_0:G_1:\ldots:G_{n-2}:G_{n-1}:G_{n},\text{ that is,}$$
$$G=G_0:\left(G_1:\ldots:\left(G_{n-2}:\left(G_{n-1}:G_{n}\right)\right)\right),$$
the values of $G_{n-1}:G_{n}$, $G_{n-2}:\left(G_{n-1}:G_{n}\right)$, $G_{n-3}\left(G_{n-2}:\left(G_{n-1}:G_{n}\right)\right)$, and so on, are iteratively calculated through equivalent positions involving the already known minimal cases. This procedure is mathematically correct, as the forms of the subordinates are irrelevant, and therefore, the known minimal forms can be used without altering the game values. Finally, we use right-to-left associativity to determine the value of the ordinal sum as a whole.

In the following sections we present a complete solution for {\sc christmas lights' fixture}, i.e., an expeditious way to compute the Grundy-value of any component. The analysis lies in finding ordinal sums like those of {\sc hackenbush}, but with carry-on moves. Theorems~\ref{case1}, \ref{case2} and \ref{case3}, in Subsection~\ref{sec:Particular}, allow us to determine the Grundy-values of three important minimal cases. Theorem~\ref{colon}, a version of the Colon Principle, is presented in Subsection~\ref{sec:Colon}. Finally, Subsection~\ref{sec:Associative} concludes the analysis with an example that showcases the theory presented here in action.

 From now on, we will use shorthands for sockets and bulbs. The notation used in this paper is as follows. 

\begin{enumerate}
  \item A sequence of $n$ consecutive sockets is designated by $|n|$.
  \item A sequence of $n$ consecutive bulbs is designated by $\overline n$.\footnote{There is a slight overlap of notation here, because we inherit the notation for an unconditional Right win as $\overline\infty$. The context is sufficiently different.}
  \item Concatenations are considered. For example $|3|\,{\overline 4}\, |2|$ designates a sequence of 3 consecutive sockets, 4 consecutive bulbs and then 2 consecutive sockets.
  \item For $0\leqslant n \leqslant\infty$, $\overline{*n}$ designates any component whose affine impartial value (explained below) equals ${*n}$ and whose leftmost piece is a bulb. One exception concerns the case $n=0$, where $\overline{*n}$ has no pieces at all.
  \item The particular case of $n=\infty$, in item 4, is highlighted as $\overline{\leftmoon}$. This component equals moon and its leftmost piece is a bulb. A minimal example is explained in Observation~\ref{obs}.
\end{enumerate}

As mentioned, the affine impartial normal play theory presented in \cite{LNS} is adequate to study rulesets with entailing moves. That general theory is used in the following sections, and can be summarized through the following list.

\begin{enumerate}
  \item[$\bullet$] \emph{Omnipresence of nimbers and moon}: Given an affine impartial game form $G$, we have a nonnegative integer $n$ such that $G=_{\Imi}*n$ or we have $G=_{\Imi}\moon$ (``$=_{\Imi}$'' is the equality of games modulo affine impartial and the moon is the game form $\{\infty\,|\,\overline{\infty}\}$; from now on, for ease, we write ``$=$'' instead of ``$=_{\Imi}$''). In the first case, we say that the {\em Grundy-value} of $G$ is $\mathcal{G}(G)=n$, and, in the second case,  we say that the {\em Grundy-value} of $G$ is $\mathcal{G}(G)=\infty$.\footnote{The meaning of this symbol should not be confused with the meaning of the same symbol in Fraenkel--Smith generalized Sprague-Grundy Theory \cite{Smith1966}. Here, we are concerned with structures with entailing moves; Fraenkel-Smith Theory  considers loopy impartial games.}
  \item [$\bullet$] \emph{Determination of the Grundy-value of $G$ from its options}: Let $G$ be an affine impartial game form, and let $\Im$ be the class of nimbers.  The set of\linebreak $G$-\emph{immediate nimbers}, denoted $\imm{G}$, is the set $S_G=\GL\cap\Im$. These are the options of $G$ that are nimbers. The set of $G$-\emph{protected nimbers}, denoted $\prot{G}$, is the set of nimbers $*n$ such that, playing first, Left wins $G+*n$ by moving to some $\{\infty\,|\,G^{L\mathcal{R}}\}+*n$ or to $\infty+*n$; although Left maybe cannot move to $*n+*n$, a winning check or a checkmate is at hand. The Grundy-value of $G$ is determined by $\mathcal{G}(G)=\mex(\mathcal{G}(S_G\cup P_G))$, where ``mex'' is the set function whose output is the minimum nonnegative integer excluded from the set.\footnote{In this paper, $\mex(\mathcal{G}(S_G\cup P_G))$ means $\mex\{\mathcal{G}(g):g\in S_G\cup P_G\}$.} Of course, if $S_G\cup P_G=\Im$, then $\mathcal{G}(G)=\infty$.  Games $G$ with options $G^L=\{\infty\,|\,{*n}\}$ are common and immediately guarantee that $\Im\setminus\{*n\}\subseteq P_G$. As these moves are carrying on to $*n$, we will use the notation $\circlearrowright^{*n}$ instead of $\{\infty\,|\,{*n}\}$. If we are making explicit Grundy-values instead of game values, we will use the notation $\circlearrowright^{n}$.

  \item [$\bullet$] \emph{Determination of the Grundy-value of a disjunctive sum, knowing the\linebreak Grundy-values of the components}: If $G$ and $H$ are affine impartial game forms, then $\mathcal{G}(G+H)=\mathcal{G}(G)\oplus \mathcal{G}(H)$, where $\oplus$ is the exclusive {\sc or} ({\sc xor}) of the binary representations of the summands if $\mathcal{G}(G)<\infty$ and $\mathcal{G}(H)<\infty$, or $\mathcal{G}(G)\oplus \mathcal{G}(H)$ results in $\infty$ if $\mathcal{G}(G)=\infty$ or $\mathcal{G}(H)=\infty$ (this operation is a natural extension of {\sc nim}-sum).
  \item [$\bullet$] \emph{Relation between the Grundy-value of $G$ and its outcome}: Given an affine impartial game form $G$, the outcome of $G$ is $\mathcal{P}$ if and only if $\mathcal{G}(G)=0$.
\end{enumerate}

Next, we present the values of all {\sc christmas lights' fixture} components with three pieces. The following section will provide closed formulas to determine the Grundy-values of these and other important particular cases.

\begin{center}
\begin{tabular}{|l|}
  \hline
  $\,$\\
  \scalebox{0.63}{\fixture{RRR}} \raisebox{0.38cm}{$=\overline{3}=\{0,*,*2\,|\,\,0,*,*2\}=*3$} \\
  \hline
    $\,$\\
  \scalebox{0.63}{\fixture{RRB}} \raisebox{0.38cm}{$=\overline{2}\, |1| =\{0,*,*2\,|\,\,0,*,*2\}=*3$} \\
  \hline
    $\,$\\
  \scalebox{0.63}{\fixture{RBR}}\raisebox{0.38cm}{$=\overline{1}\,|1|\,\overline{1}=\{0,\circlearrowright^{*2},*2\,|\,\,0,\circlearrowright^{*2},*2\}=\leftmoon$} \\
  \hline
    $\,$\\
  \scalebox{0.63}{\fixture{BRR}} \raisebox{0.38cm}{$=|1|\,\overline{2}=\{*2,*,0\,|\,\,*2,*,0\}=*3$} \\
  \hline
    $\,$\\
  \scalebox{0.63}{\fixture{BBR}} \raisebox{0.38cm}{$=|2|\,\overline{1}=\{0,0,0\,|\,\,0,0,0\}=*$} \\
  \hline
    $\,$\\
  \scalebox{0.63}{\fixture{BRB}} \raisebox{0.38cm}{$=|1|\,\overline{1}\,|1|=\{*2,*,0\,|\,\,*2,*,0\}=*3$} \\
  \hline
    $\,$\\
  \scalebox{0.63}{\fixture{RBB}} \raisebox{0.38cm}{$=\overline{1}\,|2|=\{0,*2,*2\,|\,\,0,*2,*2\}=*$} \\
  \hline
    $\,$\\
  \scalebox{0.63}{\fixture{BBB}} \raisebox{0.38cm}{$=|3|=\{0,0,0\,|\,\,0,0,0\}=*$} \\
  \hline
\end{tabular}
\end{center}

\section{Analysis of {\sc christmas lights' fixture}}\label{sec:Analysis}
A component that only has bulbs is isomorphic to a \textsc{green hackenbush string}. A component that only has sockets is a trivial {\sc she loves me she loves me not} situation, since at each move, exactly one piece is fixed. Therefore, our analysis begins with elementary positions with two types of pieces.

First, we present and prove some closed formulas that are useful for determining the Grundy-values of important particular cases. Second, we prove a kind of Colon Principle, stating that the Grundy-value of $\overline{k}\,|m|\,\overline{*i}$ only depends on the Grundy-value $i$ and not on the shape of $\overline{*i}$ ($k,m\geqslant0$, and $\infty\geqslant i\geqslant0$). Finally, we exemplify how to use right-to-left associativity to compute the Grundy-value of any component.

The following proofs are made by induction. Typically we determine $S_G$ and $P_G$ in order to compute $\mathcal{G}(G)=\mex(\mathcal{G}(S_G\cup P_G))$. Of course, the values $\mathcal{G}(S_G\cup P_G)$ are obtained through the inductive step.

\subsection{Grundy-values of $|m|\,\overline{n}$, $\overline{n}\,|m|$, and $\overline{k}\,|m|\,\overline{n}$}\label{sec:Particular}

When playing in components of the $|m|\,\overline{n}$ or $\overline{n}\,|m|$ types, sockets never appear sandwiched between two bulbs. Therefore, at all moments, the available options are {\em quiet} options,\footnote{Quiet options do not involve terminating threats, carry-on moves and so forth.} making analysis relatively simple, as illustrated in the following theorem.

\begin{theorem}\label{case1} If $m>0$, $n\geqslant 0$, $G=|m|\,\overline{n}$, and $H=\overline{n}\,|m|$, then $$\mathcal{G}(G)=\left\{\begin{array}{lll}
                                                                             n &  & \text{if } m \text{ is even} \\
                                                                             n+(-1)^n &  & \text{if } m \text{ is odd}
                                                                           \end{array}
\right.
$$
and
$$\mathcal{G}(H)=\left\{\begin{array}{lll}
                                                                             n &  & \text{if } m \text{ is even} \\
                                                                             n+1 &  & \text{if } m \text{ is odd.}
                                                                           \end{array}
\right.
$$
\end{theorem}

\begin{proof}
Let $G$ be a component of the form $|m|\,\overline{n}$. Since there are no sockets sandwiched between two bulbs, all options are quiet options. Hence, $P_G=\varnothing$, all options belong to $S_G$, and $|m|\;\overline{n}$ is a nimber.

If $m=1$ and $n=0$, we have $|1|\,\overline{0}=\{0\,|\,0\}=*$, and that is consistent with the formula. This is the base case.

Otherwise, the options of $|m|\, \overline{n}$ are $(|{m-1}|\,\overline{n}$, $|m|\,\overline{0}$, $|m|\,\overline{1}$,\ldots, $ |m|\,\overline{n-1}$.

If $m$ is even, by induction, the Grundy-values of the options are $n-1$ or $n+1$, $0$, $1$,\ldots, $n-1$. In both cases, the minimum excluded value is $n$.

If $m$ is odd and $n$ is odd, by induction, the Grundy-values of the options are $n$, $1$, $0$, $3$, $2$,\ldots, $n-4$, $n-5$, $n-2$, $n-3$, $n$. The minimum excluded value is $n-1$.

If $m$ is odd and $n$ is even, by induction, the Grundy-values of the options are $n$, $1$, $0$, $3$, $2$,\ldots, $n-3$, $n-4$, $n-1$, $n-2$. The minimum excluded value is $n+1$, and the proof is finished.\\

Let $H$ be a component of the form $\overline{n}\,|m|$. Since there are no sockets sandwiched between two bulbs, all options are quiet options. Hence, $P_H=\varnothing$, all options belong to $S_H$, and $\overline{n}\,|m|$ is a nimber.

If $m=1$ and $n=0$, we have $\overline{0}\,|1|=\{0\,|\,0\}=*$, and that is consistent with the formula. This is the base case.

Otherwise, the options of $\overline{n}\,|m|$ are $\overline{0},\ldots,\overline{n-1}$, and $\overline{n}\,|m-1|$.

If $m$ is even, by induction, the Grundy-values of the options are $0$, $1$,\ldots, $n-1$, and $n+1$. The minimum excluded value is $n$.

If $m$ is odd, by induction, the Grundy-values of the options are $0$, $1$,\ldots, $n-1$, and $n$. The minimum excluded value is $n+1$, and the proof is finished.
\end{proof}

\begin{observation}
Consider

\begin{center}\scalebox{0.63}{\fixture{BBBRR}} \raisebox{0.38cm}{$=|3|\,\overline{2}=*3$ ($m=3$ is odd and $n=2$ is even)}\end{center}

This component is trivially isomorphic to the \textsc{green hackenbush} position\footnote{Horizontal presentation.}
\begin{center}\includegraphics[width=0.25\linewidth]{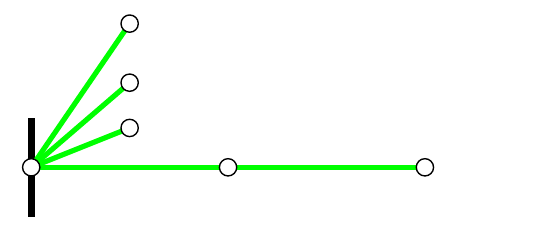}\end{center}

In fact, all components of the form $|m|\,\overline{n}$ exhibit behavior similar to \textsc{green hackenbush} positions where $m$ single edges are arranged side by side with a string of length $n$. This explains the first case of the previous theorem.\\

Consider now

\begin{center}\scalebox{0.63}{\fixture{RRBBB}} \raisebox{0.38cm}{$=\overline{2}\,|3|=*3$ ($n=2$ is even and $m=3$ is odd )}\end{center}

This component is trivially isomorphic to the \textsc{green hackenbush} position
\begin{center}\includegraphics[width=0.25\linewidth]{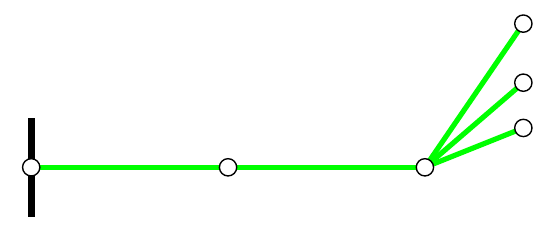}\end{center}

Indeed, components of the form $\overline{n}\,|m|$ exhibit behavior similar to \textsc{green hackenbush} positions where a string of length $n$ has $m$ single edges on the rightmost vertex. This explains the second case of the previous theorem.
\end{observation}

The following theorem, based on a minimal ``lunar situation'', already allows for the possibility of a carry-on move.

\begin{theorem}\label{case2} If $k>0$, $m>0$, and $G=\overline{k}\,|m|\,\overline{1}$, then $$\mathcal{G}(G)=\left\{\begin{array}{lll}
                                                                             \infty &  & \text{if } m=1\\
                                                                             k &  & \text{if } m>1 \text{ is odd} \\
                                                                              k+1&  & \text{if } m \text{ is even.}
                                                                           \end{array}
\right.$$
\end{theorem}

\begin{proof} If $m=1$, the socket is sandwiched between two bulbs and $G$ has one carry-on option. The options of $G$ are $\overline{0}$,\ldots, $\overline{k-1}$, a carry-on move to $\overline{k+1}$, and $\overline{k}\,|1|$. Hence, $S_{G}=\{0,*,\ldots, *(k-2),*(k-1),*(k+1)\}$, and, due to the fact that the carry-on move is $\circlearrowright^{*(k+1)}$, we have $P_{G}=\Im\setminus\{*(k+1)\}$. Therefore, $\mathcal{G}(G)=\mex(\mathcal{G}(S_{G}\cup P_{G}))=\mex(\mathcal{G}(\Im))=\infty$.

Regarding the second case, the options of $G$ are $\overline{0}$, $\overline{1}$,\ldots, $\overline{k-2}$, $\overline{k-1}$, $\overline{k}\,|m-1|\,\overline{1}$, and $\overline{k}\,|m|$. Hence, the Grundy-values of the options of $G$ are $0$, $1$,\ldots, $k-2$, $k-1$, $k+1$, and $k+1$. The penultimate is obtained by induction. The last term is obtained taking into account Theorem \ref{case1}. The minimum excluded value is $k$.

Regarding the third case, the options of $G$ are $\overline{0}$, $\overline{1}$,\ldots, $\overline{k-2}$, $\overline{k-1}$, $\overline{k}\,|m-1|\,\overline{1}$, and $\overline{k}\,|m|$. Hence, the Grundy-values of the options of $G$ are $0$, $1$,\ldots, $k-2$, $k-1$, $\infty$ or $k$, and $k$. If $m-1=1$ then the penultimate is obtained taking into account the first case of this proof; if $m-1$ is an odd integer larger than $1$ then the penultimate is obtained by induction. The last term is obtained by using Theorem \ref{case1}. The minimum excluded value is $k+1$, and the proof is finished.
\end{proof}

\begin{observation}\label{obs}
If $k=1$ and $m=1$, then we have the situation

\begin{center} \scalebox{0.63}{\fixture{RBR}} \raisebox{0.38cm}{$=\{0,\{\infty\,|\,*2\},*2\,|\,0,\{*2\,|\,\overline{\infty}\},*2\}=\leftmoon$} \end{center}
In this case, since the socket is sandwiched between two bulbs, there is a carry-on option. In fact, there is a quiet option to $*2$ and there is a carry-on option to $*2$. Consider the disjunctive sum $G+*2$; in that sum, a move to $*2+*2$ is a winning move for the first player. Consider now the disjunctive sum $G+*j$ where $j\neq 2$; in that case, the carry-on move $G+*j\rightarrow\circlearrowright^{*2}+*j$ is a winning move because the first player can continue playing on the $\mathcal{N}$-position $*2+*j$. In \cite{BerleCG1982}, the authors called this situation ``a kind of strategy stealing'' (page 406). The position $\overline{1}\,|1|\,\overline{1}$ is the simplest component whose value is equal to $\leftmoon$.\\

If $k=2$ and $m=4$ we have the situation

\begin{center}\scalebox{0.63}{\fixture{RRBBBBR}} \raisebox{0.38cm}{$=*3$}\end{center}

\vspace{-0.6cm}
\hspace{.56cm}\hspace{3.15cm}$0$\hspace{0.595cm}$1$\hspace{1.7cm}$2$\hspace{1.75cm}$2$\\

The labels on the options represent their Grundy-values. In this example, fixing the rightmost bulb makes it so that a ignorable even number of sockets remains. Therefore, in practical terms, the simplest way to approach this type of situation is to simply disregard the sockets.

If $k=2$ and $m=3$ we have the situation

\begin{center}\scalebox{0.63}{\fixture{RRBBBR}} \raisebox{0.38cm}{$=*2$}\end{center}

\vspace{-0.6cm}
\hspace{.56cm}\hspace{3.55cm}$0$\hspace{0.595cm}$1$\hspace{1.32cm}$3$\hspace{1.35cm}$3$\\

In this last example, the removals of sockets and the replacement of the bulb are reversible options. Therefore, in practice, the best approach to this type of situation is to treat the position as if it only has the bulbs on the left side of the sockets.
\end{observation}

Theorems \ref{case1} and \ref{case2} allow the analysis of components $\overline{k}\,|m|\,\overline{n}$ where $n>1$. As we will see, these cases can also be thought of as particular positions of \textsc{green hackenbush} positions.\\

\begin{theorem}\label{case3} If $k>0$, $m>0$, $n>1$, and $G=\overline{k}\,|m|\,\overline{n}$, then $$\mathcal{G}(G)=\left\{\begin{array}{lll}
                                                                              k+n&  & \text{if } m \text{ is odd}\\
                                                                              k+n+(-1)^n&  & \text{if } m \text{ is even.}
                                                                           \end{array}
\right.$$
\end{theorem}

\begin{proof} The options of $G$ are $\overline{0}$, $\overline{1}$,\ldots, $\overline{k-1}$, $\overline{k}\,|m-1|\,\overline{n}$ (if $m>1$) or $\circlearrowright^{*(k+n)}$ (if $m=1$), $\overline{k}\,|m|$, $\overline{k}\,|m|\,\overline{1}$, $\overline{k}\,|m|\,\overline{2}$,\ldots, $\overline{k}\,|m|\,\overline{n-1}$. In the following lines, the Grundy-value of the option $\overline{k}\,|m|$ is obtained by using Theorem \ref{case1} and the Grundy-values of all the other options are obtained by induction (including quiet removals of sockets).

Let $m$ be odd and $n$ be odd. If $m>1$, then we have the fundamental sets $S_{G}=\{0,\ldots, *(k-1), *(k+n-1), *(k+1), *k, *(k+2), \ldots, *(k+n-1)\}$ and $P_G=\varnothing$. If $m=1$, then, since there is one carry-on move, we have the fundamental sets $S_{G}=\{0,\ldots, *(k-1), *(k+1), *(k+2), \ldots, *(k+n-1)\}$, and $P_G=\Im\setminus\{*(k+n)\}$. In both cases, $\mathcal{G}(G)=\mex(\mathcal{G}(S_{G}\cup P_{G}))=k+n$.

Let $m$ be odd and $n$ be even. if $m>1$, then we have the fundamental sets $S_{G}=\{0,\ldots, *(k-1), *(k+n+1), *(k+1), *k, *(k+2), \ldots, *(k+n-1)\}$ and $P_G=\varnothing$. If $m=1$, then, since there is one carry-on move, we have the fundamental sets $S_{G}=\{0,\ldots, *(k-1), *(k+1), *(k+2), \ldots, *(k+n-1)\}$, and $P_G=\Im\setminus\{*(k+n)\}$. In both cases, $\mathcal{G}(G)=\mex(\mathcal{G}(S_{G}\cup P_{G}))=k+n$.

Let $m$ be even and $n$ be odd. In this case, the Grundy-values of the options of $G$ are $0$, $1$,\ldots, $k-2$, $k-1$, $k+n$ (removal of a socket), $k$, $k+1$, $k+3$, $k+2$, \ldots, $k+n-2$, $k+n-3$, and $k+n$. The minimum excluded value is $k+n-1$.

Let $m$ be even and $n$ be even. In this case, the Grundy-values of the options of $G$ are $0$, $1$,\ldots, $k-2$, $k-1$, $k+n$ (removal of a socket), $k$, $k+1$, $k+3$, $k+2$, \ldots, $k+n$, and $k+n-1$. The minimum excluded value is $k+n+1$, and the proof is finished.\\
\end{proof}

\begin{observation}
Consider the component where $k=2$, $m=1$, and $n=4$, that is,

\begin{center}\scalebox{0.63}{\fixture{RRBRRRR}} \raisebox{0.38cm}{$\,=*6$} \end{center}

\vspace{-0.6cm}
\hspace{.56cm}\hspace{3.1cm}$0$\hspace{0.595cm}$1$\hspace{0.55cm}$\circlearrowright^{6}$\hspace{0.39cm}$3$\hspace{0.59cm}$\leftmoon$\hspace{0.55cm}$4$\hspace{0.58cm}$5$\\

In this case, the Grundy-value of the component coincides with the total number of bulbs. This happens whenever $m=1$ and $n>1$. Consequently, in terms of game practice, the simplest approach is to make the carry-on move and, if appropriate, play again in the same component.

\begin{center}\scalebox{0.63}{\fixture{RRRRRR}} \raisebox{0.38cm}{$\,=*6$} \end{center}

\vspace{-0.6cm}
\hspace{.56cm}\hspace{3.5cm}$0$\hspace{0.595cm}$1$\hspace{0.595cm}$2$\hspace{0.59cm}$3$\hspace{0.55cm}$4$\hspace{0.58cm}$5$\\

This can be thought of as a {\sc green hackenbush} string with a ``ghost edge''.

\begin{center}\includegraphics[width=0.9\linewidth]{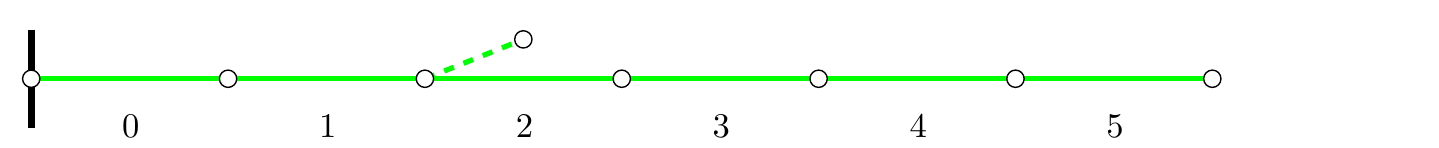}\end{center} In practice, the ghost edge does not exist,\footnote{Ghosts do not exist.} since it corresponds to a carry-on move that keeps the right to play to the player who makes it. Consider now the position where $k=2$, $m=2$, and $n=5$, that is,

\begin{center}\scalebox{0.63}{\fixture{RRBBRRRRR}} \raisebox{0.38cm}{$\,=*6$} \end{center}

\vspace{-0.6cm}
\hspace{.56cm}\hspace{2.35cm}$0$\hspace{0.595cm}$1$\hspace{0.95cm}$7$\hspace{0.98cm}$2$\hspace{0.58cm}$3$\hspace{0.57cm}$5$\hspace{0.57cm}$4$\hspace{0.57cm}$7$\\

This example is more difficult to analyze. Applying the previous theorem, given that $n$ is odd, we obtain the Grundy-value $2+5-1=6$. However, due Theorem~\ref{case2}, there is a ``perturbation'' related to the two bulbs following the sockets, so it is no longer easy to have an intuition about the Grundy-values of the options. To overcome this problem, we suggest the mnemonic \emph{Double Jump} which consists of establishing again a link with a {\sc green hackenbush} position, but placing the second edge two vertices to the right what would be
expected.

\begin{center}\includegraphics[width=0.9\linewidth]{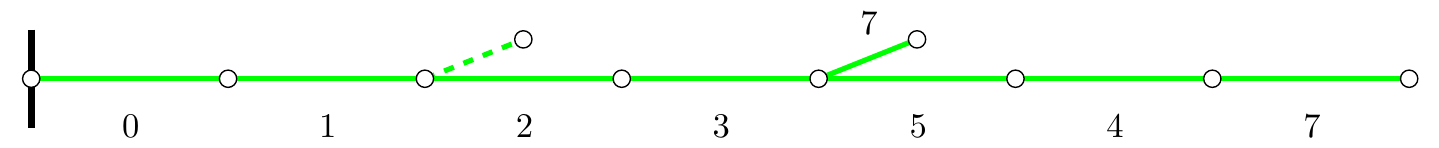}\end{center}

Algebraically speaking, this action regularizes the perturbation. Then, it is only necessary to determine the Grundy-values of the options of that {\sc green hackenbush} position, always keeping in mind that the ghost does not exist.\\

For larger values of $m$, the mnemonic still works. However, starting from the second one, ``non-ghost'' edges are placed on the expected vertex. For example, consider $k=2$, $m=3$, and $n=5$, that is,

\begin{center}\scalebox{0.63}{\fixture{RRBBBRRRRR}} \raisebox{0.38cm}{$\,=*7$} \end{center}

\vspace{-0.6cm}
\hspace{.56cm}\hspace{1.97cm}$0$\hspace{0.595cm}$1$\hspace{1.35cm}$6$\hspace{1.35cm}$3$\hspace{0.58cm}$2$\hspace{0.57cm}$4$\hspace{0.57cm}$5$\hspace{0.57cm}$6$\\

In that case, the related {\sc green hackenbush} position is
\begin{center}\includegraphics[width=0.9\linewidth]{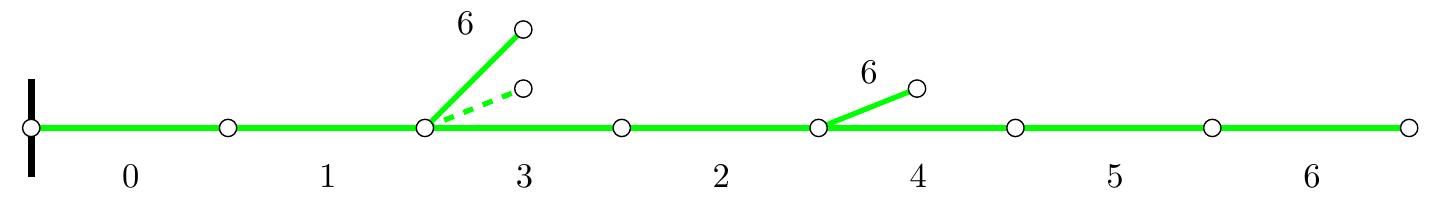}\end{center}
and, for $\overline{2}\,|4|\,\overline{5}$ and $\overline{2}\,|5|\,\overline{5}$, the related {\sc green hackenbush} positions are
\begin{center}\includegraphics[width=0.9\linewidth]{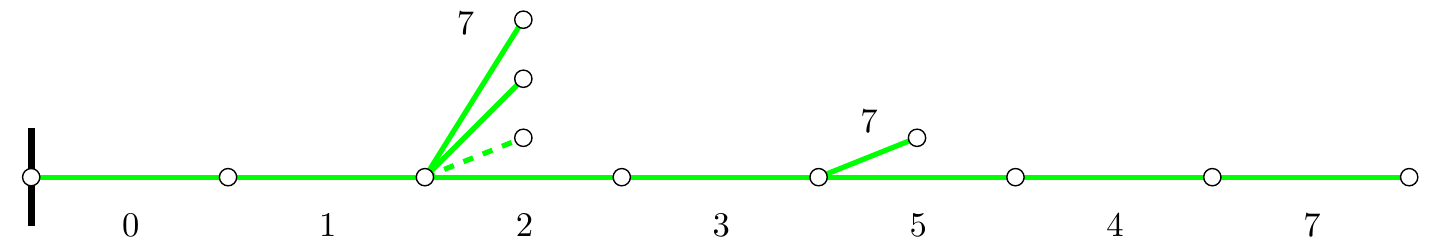}\end{center}and\begin{center}\includegraphics[width=0.9\linewidth]{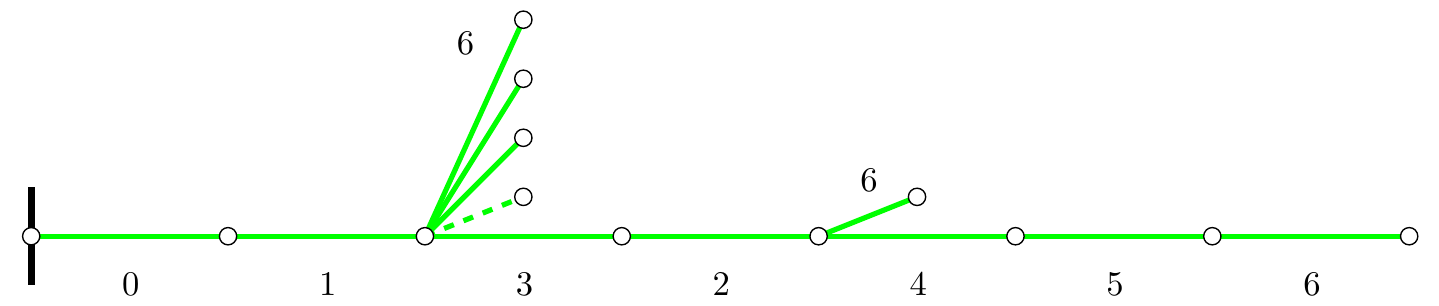}\raisebox{0.07cm}{$\,$.}\end{center}
\end{observation}

\vspace{0.3cm}
\subsection{Moonlight Theorem and Colon Principle}\label{sec:Colon}

The initial results of this section are related to the Colon Principle, but they only apply to components with a finite Grundy-value. Nevertheless, once we establish the Moonlight Theorem, we prove that the Colon Principle is applicable to all cases.

\begin{lemma}\label{collemma2} If $i<s$ are two integers, we have the following:
\begin{enumerate}
  \item All carry-on moves on $\overline{*i}$ are $\circlearrowright^{\overline{*i}'}$ where $\overline{*i}'$ is identical to the fixture $\overline{*i}$ in every way, except for the absence of the socket removed with the carry-on move;
  \item In the game $\overline{*i}+\overline{*s}$, the first player can win by making all possible carry-on moves on the second component, followed by a move to $\overline{*i}+\overline{*i}'$ in that same component.
\end{enumerate}
\end{lemma}

\begin{proof}
Starting with the first item, note that a carry-on move on $\overline{*i}$ cannot be $\circlearrowright^{\overline{\leftmoon}}$, otherwise the Grundy-value of $\overline{*i}$ would be $\infty$ instead of $i$. On the other hand, a carry-on move on $\overline{*i}$ also cannot be $\circlearrowright^{\overline{*j}}$ with $j\neq i$, otherwise we would have $\Im\setminus\{*j\}\subseteq P_{\overline{*i}}$, and $\mex(\mathcal{G}(S_{\overline{*i}}\cup P_{\overline{*i}}))$ could be $j$ or $\infty$, but not $i$. Therefore, a carry-on move on $\overline{*i}$ can only be $\circlearrowright^{\overline{*i}'}$.

The second item is a consequence of the first. That is, after making all possible carry-on moves on the second component, a position $\overline{*i}+\overline{*s}'$ is obtained in which the second component no longer has any carry-on moves. Since its Grundy-value is still $s$, at that point, the second component must have a quiet move to $\overline{*i}'$. Naturally, the move to $\overline{*i}+\overline{*i}'$ is a winning move since its Grundy-value is $i\oplus i=0$.
\end{proof}

\begin{lemma}\label{collemma} Let $k\geqslant0$, $m\geqslant0$, and $i<\infty$. If $\overline{k}\,|m|\,\overline{*i}$ is not the moon, then the Grundy-value of $\overline{k}\,|m|\,\overline{*i}$ does not depend on the shape of $\overline{*i}$.
\end{lemma}

\begin{proof}
Let $\overline{*i}$ and $\overline{*i}'$ be two components of different shapes, but with the same Grundy-value $i<\infty$. First, suppose that $m=0$. In this case, $\overline{k}\,\overline{*i}+\overline{k}\,\overline{*i}'\in\mathcal{P}$ because the second player can play on the rightmost parts as they were playing the $\mathcal{P}$-position $\overline{*i}+\overline{*i}'$. Hence, $\overline{k}\,\overline{*i}=\overline{k}\,\overline{*i}'$.

Now, let us prove that we cannot have $k>0$ and $m=i=1$. In $\overline{k}\,|1|\,\overline{*}$, by Theorem \ref{case1}, a player would be able to make a move to $\overline{k}\,|1|=*(k+1)$. On the other hand, in $\overline{k}\,|1|\,\overline{*}$, a player could also make a carry-on move to $\overline{k}\,\overline{*}$. However, we have already shown in the previous paragraph that $\overline{k}\,\overline{*}=\overline{k}\,\overline{1}$. Since $\overline{k}\,\overline{1}=*(k+1)$, that carry-on move would lead to $\overline{k}\,\overline{*}=*(k+1)$. Consequently, we have $*(k+1)\in S_{\overline{k}\,|1|\,\overline{*}}$, but at the same time, $\Im\setminus\{*(k+1)\}\subseteq P_{\overline{k}\,|1|\,\overline{*}}$. Together, these two facts imply that $\mathcal{G}(\overline{k}\,|1|\,\overline{*})=\mex(\mathcal{G}(S_{\overline{k}\,|1|\,\overline{*}}\cup P_{\overline{k}\,|1|\,\overline{*}}))=\mex(\mathcal{G}(\Im))=\infty$, therefore, it follows that $\overline{k}\,|1|\,\overline{*}=\leftmoon$, contradicting our assumption.

For the general case, let us prove that $\overline{k}\,|m|\,\overline{*i}+\overline{k}\,|m|\,\overline{*i}'\in\mathcal{P}$. Essentially, the strategy of the second player is to play on the rightmost parts as if they were playing the $\mathcal{P}$-position $\overline{*i}+\overline{*i}'$.

If the first player makes a quiet move on one of the rightmost parts $\overline{*i}$ or $\overline{*i}'$, then the second player responds also on the rightmost parts with their winning line of $\overline{*i}+\overline{*i}'$. As a result, a position like $\overline{k}\,|m|\,\overline{*j}+\overline{k}\,|m|\,\overline{*j}'$ is achieved. By induction on the number of pieces, that position is a $\mathcal{P}$-position and the second player wins.

If the first player makes a carry-on move on one of the rightmost parts, say $\overline{*i}$, then, by Lemma \ref{collemma2}, that move carries on to some $\overline{*i}''$. Hence, the first player has the turn again in the position $\overline{k}\,|m|\,\overline{*i}''+\overline{k}\,|m|\,\overline{*i}'$. Since by induction on the number of pieces that position is a $\mathcal{P}$-position, the second player wins.

If the first player replaces a bulb from one of the leftmost parts $\overline{k}$ of one of the components, the second player mimics that move on the other component, obtaining a position $\overline{k-j}+\overline{k-j}=*(k-j)+*(k-j)\in\mathcal{P}$.

If $m>1$ and the first player removes one of the $m$ sockets from one of the components, the second player mimics that move on the other component, obtaining the position $\overline{k}\,|m-1|\,\overline{*i}+\overline{k}\,|m-1|\,\overline{*i}'$. Once more, by induction on the number of pieces,  $\overline{k}\,|m-1|\,\overline{*i}+\overline{k}\,|m-1|\,\overline{*i}'$ is a $\mathcal{P}$-position and the second player wins.

If $k>0$, $m=1$ and the first player makes a carry-on move by removing that single socket, followed by a sequence of moves that allows the second player to obtain a position such as  $\overline{k}\,\overline{*j}+\overline{k}\,\,\overline{*j}'$ or a position such as $\overline{k-j}+\overline{k-j}$, the second player also wins. The first case is a $\mathcal{P}$-position, as explained in the first paragraph of this proof.  In the second case, we have $\overline{k-j}+\overline{k-j}=*(k-j)+*(k-j)\in\mathcal{P}$. The only scenario where the second player is unable to reach such positions is when $k>0$, $m=1$, and the first player makes a carry-on move on one component to $\overline{k}\,\overline{*i}+\overline{k}\,|1|\,\overline{*i}'$, followed by a move on the second component to $\overline{k}\,\overline{*i}+\overline{k}\,|1|$.  This is the only case where the second player cannot use the winning strategy of $\overline{*i}+\overline{*i}'$ (moving from $\overline{*i}+0$ to $0+0$) because they no longer have access to a preliminary carry-on move on the second component before doing it. However, even in this case, since we already know that $i\neq 1$, the second player can still win by playing on the rightmost parts by choosing the move that corresponds to the winning move of $\overline{*i}+|1|$.
\end{proof}

The following theorem states that when a moon appears as a part of a component, the value of the component as a whole is also equal to the moon. In other words, regarding this type of concatenations, the moon is an absorbing element.

\begin{theorem}[Moonlight Theorem]\label{moonlight} If $k\geqslant0$ and $m\geqslant0$ then $\overline{k}\,|m|\,\overline{\leftmoon}=\leftmoon$.\footnote{When we are exposed to the moonlight, regardless of how far away the moon is, it is impossible to ignore its presence.}
\end{theorem}

\begin{proof}
Suppose there exists a component $\overline{\leftmoon}$ such that $\mathcal{G}(\overline{k}\,|m|\,\overline{\leftmoon})<\infty$. Additionally, assume that $\overline{\leftmoon}$ is composed of the minimal number of pieces possible. The component $\overline{\leftmoon}$ must include carry-on moves, otherwise its game value would not be the moon.

A carry-on move of $\overline{\leftmoon}$ cannot be $\{\infty\,|\,\overline{\leftmoon}'\}$, where $\overline{\leftmoon}'$ is identical to $\overline{\leftmoon}$ in every way, except for the absence of the socket removed with the carry-on move. If it were, then, due to the minimal assumption, $\overline{k}\,|m|\,\overline{\leftmoon}'$ would be the moon. And, because of that, $\overline{k}\,|m|\,\overline{\leftmoon}$ would have a carry-on move to the moon, contradicting the assumption that $\mathcal{G}(\overline{k}\,|m|\,\overline{\leftmoon})<\infty$.

Suppose now that all carry-on moves of $\overline{\leftmoon}$ have the form $\{\infty\,|\,\overline{*i}\}$ for some $i<\infty$. Because the game value of $\overline{\leftmoon}$ is the moon, a player must also have a move to $\overline{*i}'$. Therefore, in $\overline{k}\,|m|\,\overline{\leftmoon}$, a player has a move to $\overline{k}\,|m|\,\overline{*i}'$ and a carry-on move to $\overline{k}\,|m|\,\overline{*i}$. If $\overline{k}\,|m|\,\overline{*i}=\leftmoon$, then $\overline{k}\,|m|\,\overline{\leftmoon}$ is the moon, contradicting $\mathcal{G}(\overline{k}\,|m|\,\overline{\leftmoon})<\infty$. On the other hand, if $\overline{k}\,|m|\,\overline{*i}=*w$ with $w<\infty$, by Lemma \ref{collemma}, that fact does not depend on the shape of $\overline{*i}$ and $\overline{k}\,|m|\,\overline{*i}'$ is also equal to $*w$. In this way, we have $*w\in S_{\overline{k}\,|m|\,\overline{\leftmoon}}$, $\Im\setminus\{*w\}\subseteq P_{\overline{k}\,|m|\,\overline{\leftmoon}}$, and $\mex(\mathcal{G}(S_{\overline{k}\,|m|\,\overline{\leftmoon}}\cup P_{\overline{k}\,|m|\,\overline{\leftmoon}}))=\infty$. Thus, $\overline{k}\,|m|\,\overline{\leftmoon}$ is the moon, contradicting once again the assumption $\mathcal{G}(\overline{k}\,|m|\,\overline{\leftmoon})<\infty$.

Finally, suppose that in $\overline{\leftmoon}$ there are at least two carry-on moves $\{\infty\,|\,\overline{*i}\}$ and $\{\infty\,|\,\overline{*s}\}$ with $i<s<\infty$. If so, in $\overline{k}\,|m|\,\overline{\leftmoon}$ a player has a carry-on move to $\overline{k}\,|m|\,\overline{*i}$ and a carry-on move to $\overline{k}\,|m|\,\overline{*s}$. Neither $\overline{k}\,|m|\,\overline{*i}$ nor $\overline{k}\,|m|\,\overline{*s}$ is the moon, or else $\overline{k}\,|m|\,\overline{\leftmoon}$ would be the moon. Also, note that $\overline{k}\,|m|\,\overline{*i}$ is not equal to $\overline{k}\,|m|\,\overline{*s}$ since the first player wins $\overline{k}\,|m|\,\overline{*i}+\overline{k}\,|m|\,\overline{*s}$. Since $i<s<\infty$, the first player can force $\overline{k}\,|m|\,\overline{*i}+\overline{k}\,|m|\,\overline{*i}'$ with the winning strategy of $\overline{*i}+\overline{*s}$ given by Lemma \ref{collemma2}, winning the game because $\overline{k}\,|m|\,\overline{*i}+\overline{k}\,|m|\,\overline{*i}'$ is a $\mathcal{P}$-position by Lemma \ref{collemma}. Thus, a player has carry-on moves to two distinct nimbers. This fact implies that $P_{\overline{k}\,|m|\,\overline{\leftmoon}}=\Im$ and $\overline{k}\,|m|\,\overline{\leftmoon}$ is the moon, definitively contradicting the assumption $\mathcal{G}(\overline{k}\,|m|\,\overline{\leftmoon})<\infty$.
\end{proof}

Now, we are ready to establish the Colon Principle in general terms.

\begin{theorem}[Colon Principle]\label{colon} Let $k\geqslant0$, $m\geqslant0$, and $\infty\geqslant i\geqslant0$. Then the Grundy-value of $\overline{k}\,|m|\,\overline{*i}$ does not depend on the shape of $\overline{*i}$.
\end{theorem}

\begin{proof}
Suppose first that $i$ is finite. If $\overline{k}\,|m|\,\overline{*i}$ is not equal to the moon, then Lemma~\ref{collemma} guarantees that the shape of $\overline{*i}$ is irrelevant. On the other hand, if $\overline{k}\,|m|\,\overline{*i}$ is equal to the moon, and if there is another component $\overline{*i}'$ such that the Grundy-value of $\overline{k}\,|m|\,\overline{*i}'$ is finite, then we have a contradiction with Lemma~\ref{collemma}. Therefore, such $\overline{*i}'$ cannot exist, and the shape of $\overline{*i}$ is again irrelevant.
Suppose now that $i$ is infinite, i.e., $\overline{*i}$ is equal to the moon. In that case, by Theorem~\ref{moonlight}, $\overline{k}\,|m|\,\overline{*i}$ is equal to the moon and, once more, the shape of $\overline{*i}$ is irrelevant.
\end{proof}

\subsection{Use of right-to-left associativity}\label{sec:Associative}

Let $C_0$,\ldots, $C_n$ be pieces such that $C_i=\overline{k}$ or $\mathbold{C_i}=|k|$. To determine the Grundy-value of the component $C_0\,\mathbold{C_1}\,\ldots\,\mathbold{C_{n-1}}\,C_{n}$, we can use right-to-left associativity:
$$C_0\,\mathbold{C_{1}}\,\left(C_2\,\ldots\,\left(C_{n-4}\,\mathbold{C_{n-3}}\,\left(C_{n-2}\,\mathbold{C_{n-1}}\,C_{n}\right)\right)\right).$$

We begin the computation by applying Theorem \ref{case1}, Theorem \ref{case2}, or Theorem \ref{case3} to $C_{n-2}\,\mathbold{C_{n-1}}\,C_{n}$. If this is the moon, then the entire component is the moon (Moonlight Theorem). If $C_{n-2}\,\mathbold{C_{n-1}}\,C_{n}=*j$ with $j<\infty$, the Colon Principle allows us to replace $C_{n-2}\,\mathbold{C_{n-1}}\,C_{n}$ with $\overline{j}$. After that replacement, it is possible to apply the theorems again to compute $C_{n-4}\,\mathbold{C_{n-3}}\,\left(C_{n-2}\,\mathbold{C_{n-1}}\,C_{n}\right)$ as it was $C_{n-4}\,\mathbold{C_{n-3}}\,\overline{j}$. And so on, until reaching the leftmost piece of the component. As an example, consider the following exercise involving the disjunctive sum $\overline{2}\,|3|\,\overline{5}\,|1|\,\overline{4}\,|2|\,\overline{2}\,+\,\overline{2}\,|1|\,\overline{3}\,|1|\,\overline{1}\,|3|\,\overline{1}$.\\

\begin{figure}[htb!]
\caption{The electricians are repairing a fixture that is in terrible condition, so the game is going to be interesting! Who wins, the $\mathcal{P}$revious player or the $\mathcal{N}$ext player? If it is the $\mathcal{N}$ext player, how?}
\label{fig:*5}

\vspace{0.5cm}
\noindent
\scalebox{0.335}{\fixture{PYCRRBBBRRRRRBRRRRBBRRSRRBRRRBRBBBR}}
\end{figure}

\vspace{0.5cm}
On one hand, the disjoint component on the left is $\overline{2}\,|3|\,\overline{5}\,|1|\,\overline{4}\,|2|\,\overline{2}$. Hence, by using right-to-left associativity, and applying Theorem \ref{case2}, Theorem \ref{case3}, and the Colon Principle, we have

$$\underbrace{(\overline{2}\,|3|\,\underbrace{(\overline{5}\,|1|\,\underbrace{(\overline{4}\,|2|\,\overline{2})}_{*7})}_{*12})}_{*14}$$\\

On the other hand, the disjoint component on the right is $\overline{2}\,|1|\,\overline{3}\,|1|\,\overline{1}\,|3|\,\overline{1}$. Hence, by using right-to-left associativity, and applying Theorem \ref{case2}, Theorem \ref{case3}, the Colon Principle, and the Moonlight Theorem, we have

$$\underbrace{(\overline{2}\,|1|\,\underbrace{(\overline{3}\,|1|\,\underbrace{(\overline{1}\,|3|\,\overline{1})}_{*})}_{\leftmoon})}_{\leftmoon}$$\\

Note that, since $\overline{3}\,|1|\,\overline{1}\,|3|\,\overline{1}=\leftmoon$, we immediately know that the entire component is equal to $\leftmoon$.

Therefore, the disjunctive sum is an $\mathcal{N}$-position, and its game value is equal to $*14+\leftmoon=\leftmoon$. The first player can win by playing a sequence of carry-on moves. She can start with two carry-on moves on the right component, moving to $\overline{2}\,|3|\,\overline{5}\,|1|\,\overline{4}\,|2|\,\overline{2}+\overline{6}\,|3|\,\overline{1}$, which is equal to $*14+*6$. After that, she can carry on to $\overline{2}\,|3|\,\overline{9}\,|2|\,\overline{2}+\overline{6}\,|3|\,\overline{1}$ on the left component, maintaining the sum in $*14+*6$. Finally, she can play the quiet winning move on the left component to $\overline{2}\,|3|\,\overline{4}+\overline{6}\,|3|\,\overline{1}=*6+*6=0$.\\

After the game is over, the electricians can admire their work illuminated with $23$ shining bulbs!\\

\vspace{0.2cm}
\noindent
\scalebox{0.49}{\fixture{PYCGOYSCGOYSCGOYSCGOYSCG}}

\end{document}